\documentclass[11pt, reqno]{amsart}
\usepackage{leading}
\usepackage{blindtext}
\usepackage{bm}
\usepackage{lineno}
\newcommand{\todo}[1]{}

\usepackage[new]{old-arrows}

\makeatletter
\DeclareMathAlphabet{\mathpzc}{OT1}{pzc}{m}{it}
\usepackage[margin=2.3cm]{geometry}
\usepackage{amsmath, amssymb, amsfonts, amstext, verbatim, amsthm, mathrsfs}
\usepackage[mathcal]{eucal}
\usepackage{microtype}
\usepackage[colorlinks=true,linkcolor=blue,citecolor=red,urlcolor=blue,citebordercolor={0 0 1},urlbordercolor={0 0 1},linkbordercolor={0 0 1}]{hyperref} 
\usepackage{comment}
\DeclareSymbolFont{bbold}{U}{bbold}{m}{n}
\DeclareSymbolFontAlphabet{\mathbbold}{bbold}

\usepackage[usenames]{color}
\usepackage{aliascnt} 
\usepackage{enumerate}
\usepackage{xspace}
\usepackage{tikz-cd}
\usepackage[nameinlink]{cleveref}

\theoremstyle{plain}

\newcommand{\refnewtheoremn}[4]{
	\newaliascnt{#1}{#2}
	\newtheorem{#1}[#1]{#3}
	\aliascntresetthe{#1}
	\expandafter\providecommand\csname #1autorefname\endcsname{#4}}

\newcommand{\refnewtheorem}[3]{\refnewtheoremn{#1}{#2}{#3}{#3}}

\def\makeCal#1{
	\expandafter\newcommand\csname c#1\endcsname{\mathcal{#1}}}
\def\makeBB#1{
	\expandafter\newcommand\csname b#1\endcsname{\mathbb{#1}}}
\def\makeFrak#1{
	\expandafter\newcommand\csname f#1\endcsname{\mathfrak{#1}}}

\count@=0
\loop
\advance\count@ 1
\edef\y{\@Alph\count@}
\expandafter\makeCal\y
\expandafter\makeBB\y
\expandafter\makeFrak\y
\ifnum\count@<26
\repeat
\newtheorem{theorem}{Theorem}[section]

\refnewtheorem{lemma}{thm}{Lemma}
\refnewtheorem{assumptions}{thm}{Assumptions}
\refnewtheorem{examples}{thm}{Examples}
\refnewtheorem{claim}{thm}{Claim}
\refnewtheorem{cor}{thm}{Corollary}
\refnewtheorem{conj}{thm}{Conjecture}
\refnewtheorem{prop}{thm}{Proposition}
\refnewtheorem{proposition}{thm}{Proposition}
\refnewtheorem{quest}{thm}{Question}
\refnewtheorem{assumption}{thm}{Assumption}

\theoremstyle{definition}
\refnewtheorem{remark}{thm}{Remark}
\refnewtheorem{remarks}{thm}{Remarks}
\refnewtheorem{defn}{thm}{Definition}
\refnewtheorem{definition}{thm}{Definition}
\refnewtheorem{notn}{thm}{Notation}
\refnewtheorem{const}{thm}{Construction}
\refnewtheorem{example}{thm}{Example}
\refnewtheorem{fact}{thm}{Fact}
\refnewtheorem{problem}{thm}{Problem}
\refnewtheorem{convention}{thm}{Convention}

\newcommand{\spec}{\operatorname{spec}}
\newcommand{\one}{\boldsymbol 1}
\newcommand{\avec}{\boldsymbol a}

\newcommand{\tagref}[1]{\href{https://stacks.math.columbia.edu/tag/#1}{Tag #1}}

\newcommand{\ul}{Ulrich }

\newcommand{\wrt}{with respect to }

\newcommand{\rk}{\operatorname{rk}}

\numberwithin{equation}{subsection}
\usepackage{expl3}
\usepackage[english]{babel}
\usepackage{xparse}
\ExplSyntaxOn 

\NewDocumentCommand\alotofdots{ D(){\ldotp} O{3}}
{\mathinner{\prg_replicate:nn{#2}{#1}}}
\ExplSyntaxOff

\begin{document}
	
	\title[Existence of Ulrich and weakly Ulrich bundles in positive characteristic]{Existence of Ulrich bundles on polarized surfaces and weakly Ulrich bundles on polarized varieties in positive characteristic}
	
	\author{Jagadish Pine}
	
	\address{Indian Institute of Science Education and Research, Pune,  Main Academic, Dr Homi Bhabha Rd, Pashan, Maharashtra 411008}
	
	\email{math.jagadishpine@gmail.com, jagadish.pine@iiserpune.ac.in}

	\date{}
	
	\keywords{Ulrich bundles; Weakly Ulrich bundles; Eisenbud-Schreyer existence problem; Positive characteristic; Frobenius morphism; Chow forms}
	
	\subjclass[2020]{ 14G17 · 13A35 · 14J60 }

	\begin{abstract}
		In their seminal paper \cite{eisenbud2003resultants}, Eisenbud and Schreyer asked whether every smooth projective variety admits an Ulrich bundle with respect to each very ample polarization. This question has stimulated extensive research for more than two decades.
		Recently, Anghel \cite{anghel2026smooth} constructed counterexamples over $\mathbb C$ in dimension two, and Lopez \cite{lopez2026existence} extended this nonexistence result to higher dimensions.
		In this article, we study the existence problem for sufficiently large multiples of a fixed polarization and derive applications in positive characteristic. Over an algebraically closed field of positive characteristic, we prove that every polarized smooth surface admits an Ulrich bundle and every polarized smooth variety admits a weakly Ulrich bundle. Consequently, the Chow forms of surfaces (resp. varieties) admit linear determinantal (resp. determinantal) representations. Moreover, for every such polarized surface, the cone of cohomology tables of coherent sheaves coincides with that on $(\mathbb{P}^2, \mathcal{O}_{\mathbb{P}^2}(1))$.
	\end{abstract}

\maketitle
	
\section{Introduction}
Ulrich and weakly Ulrich bundles were introduced into algebraic geometry
by Eisenbud and Schreyer \cite{eisenbud2003resultants} in their work on
explicit formulas for resultants of homogeneous forms and Chow forms
of projective varieties. A vector bundle $E$ on an $n$-dimensional smooth
projective variety $X$ is Ulrich with respect to a very ample divisor $H$
if $H^\bullet(X,E(-tH))=0$ for $t=1,\ldots,n$. A useful generalization is the notion of a weakly Ulrich bundle, defined
by weaker cohomological vanishing conditions; see
Definition \ref{defweak}.

In the same paper, Eisenbud and Schreyer
\cite[p.~543]{eisenbud2003resultants} asked whether every smooth
projective variety admits an Ulrich bundle with respect to any
given very ample polarization. This question, which they described
as ``of greatest importance'' among the problems they left open,
has motivated extensive research.  Affirmative answers have been obtained for several classes of
polarized varieties.

However, in a recent preprint \cite{anghel2026smooth}, Anghel constructed
smooth complex projective surfaces with very ample polarization
admitting no Ulrich bundle. Building on these surfaces, Lopez obtained counterexamples in every
dimension at least $3$ \cite{lopez2026existence}.

In this article, we prove the existence of Ulrich bundles on every smooth polarized projective
surface and of weakly Ulrich bundles on every smooth polarized projective variety defined over algebraically closed fields of positive characteristic.

Ulrich bundles are slope semistable with respect to the polarization \cite[Theorem 2.9(a)]{casanellas2012stable}.
Consequently, in characteristic zero, an Ulrich bundle $E$ of rank $r$
on a smooth projective surface satisfies the Bogomolov-Gieseker inequality $\Delta(E):=2r c_2(E)-(r-1)c_1(E)^2\geq 0$.
This inequality is a crucial ingredient in Anghel's nonexistence argument, but need not be true for Ulrich bundles in positive characteristic.
In positive characteristic, it holds for strongly semistable bundles
but may fail for semistable bundles; see
\cite[p. 917]{langer2015bogomolov}.
Moreover, Ulrich bundles need not be strongly semistable.
For example, on a smooth projective ordinary curve $C$, the sheaf of exact $1$-forms $B_C^1(H)$ is Ulrich for any very ample divisor $H$ \cite[Corollary 2.2(1)]{MR3920911}, but $B_C^1(H)$ fails to be strongly semistable when genus
at least $2$ and characteristic $p>2$ \cite[Remark~4.1.2(2)]{raynaud1982sections}. Similar examples exist on $C \times C$.

\begin{theorem}\label{thm1}
	Let $S$ be a smooth projective surface over an algebraically closed
	field of characteristic $p>0$, and let $H$ be a very ample divisor on $S$.
	There exists an integer $e_0\geq 1$ such that, for every integer
	$e\geq e_0$, the polarized surface $(S,H)$ admits an Ulrich bundle
	of rank $2p^{2e}$.
\end{theorem}

We construct Ulrich bundles on $(S,H)$ using Frobenius direct images. The key step is an asymptotic existence theorem for rank-two Ulrich bundles on $(S,qH)$, which we establish in positive characteristic for all sufficiently large integers $q$. This asymptotic theorem extends the result of Coskun and Huizenga \cite[Theorem 4.3]{CH20} for smooth complex projective surfaces. Using the above theorem, in Section \ref{productcons} we provide examples of polarized varieties admitting Ulrich bundles, including products with surfaces and certain projective bundles over surfaces.

Joshi studied the ACM and Ulrich
properties of the bundle $B_X^1$, and also Frobenius direct images of
the structure sheaf in \cite{MR3920911, Jos21}.
He constructed weakly Ulrich
bundles on every smooth projective Frobenius split surface and
threefold \cite[Theorem 5.4]{Jos21}. He also obtained examples on
certain surfaces of general type, and certain surfaces in $\bP^3$ in \cite[Theorem 3.1, Theorem 4.1]{Jos21}. 

Our result on weakly Ulrich bundles establishes their existence in arbitrary dimension.

\begin{theorem}\label{thm2}
	Let $X$ be a smooth projective variety of dimension $n$ over an algebraically closed field $k$ of characteristic $p>0$, and let $H$ be very ample. There is an integer $e_0$ such that, for every $e\geq e_0$, the pair $(X,H)$ carries a weakly Ulrich vector bundle $E_e$ of rank $n!p^{en}$.
	\end{theorem}

Theorem \ref{thm2} follows, via Frobenius direct image, from
the following asymptotic existence result, which holds in
arbitrary characteristic.

\begin{theorem}\label{thm3}
	Let $X$ be a smooth projective variety of dimension $n$ over an algebraically closed field of arbitrary characteristic, and let $H$ be very ample. There is an integer $q_0$ such that $(X,qH)$ carries a weakly Ulrich vector bundle $B_q$ of rank $n!$ for every integer $q\geq q_0$.
\end{theorem}

We outline the construction. Fix a finite linear projection
$\pi:X\longrightarrow\mathbb P^n$ and a finite flat map
$\rho:(\mathbb P^1)^n\longrightarrow\mathbb P^n$ such that $\rho^*\mathcal O_{\mathbb P^n}(1) \cong \mathcal O_{\mathbb P^1}(1)\boxtimes\cdots\boxtimes \mathcal O_{\mathbb P^1}(1)$.
For a suitable family of line bundles $L_q$ on $(\mathbb P^1)^n$, set
\[
B_q=\pi^*(\rho_*L_q).
\]
We establish the cohomological vanishings that make $B_q$
weakly Ulrich with respect to $qH$ for every sufficiently large
integer $q$, in arbitrary characteristic.

The choice of $L_q$ arise from construction of Ulrich bundles on product of varieties;
see Section \ref{productcons} and
\cite[Proposition 3.1]{beauville2018introduction}.
The underlying construction also appears in the work of Iyengar, Ma
and Walker on lim Ulrich sequences \cite[\S 7]{IMW22};
see also \cite{IMW23}.
The weakly Ulrich property established here does not follow from the
lim Ulrich condition: a lim Ulrich sequence may have no weakly Ulrich
terms, and a sequence of weakly Ulrich bundles need not be lim Ulrich.
We explain this distinction in Section \ref{weaklimrel},
with Examples \ref{limweakex1} and \ref{limweakex2}.

Taking $q=p^e$ for $e$ sufficiently large, the Frobenius direct image
$E_e=F^e_*B_q$ gives the desired weakly Ulrich bundle on $(X,H)$.
We also show that the constructed weakly Ulrich bundles are initialized;
see Proposition \ref{initialized}.
The bundles $B_q$ and $E_e$ need not be Ulrich with respect to
$qH$ and $H$, respectively; see Remark \ref{neednotul}.

Joshi introduced the notion of almost Ulrich bundles, which are both ACM and weakly Ulrich (see Definition \ref{almostul}), and proved their existence on Frobenius split threefolds \cite[Theorem 5.4]{Jos21}. The weakly Ulrich bundles $E_e$ need not be almost Ulrich; see Remark \ref{nonalmostul}.

 \begin{theorem} \label{thm4}
Every smooth projective threefold over an algebraically closed field of positive characteristic carries an almost Ulrich vector bundle \wrt every very ample divisor $H$.
\end{theorem}

The key step in proving the existence of almost Ulrich bundles is the construction of an acyclic bundle $B$ on every projective variety $X$ in Theorem \ref{acyclic}. When $X$ is smooth of dimension $3$ and defined over an algebraically closed field of positive characteristic, a suitable twist of $F^e_*B$ by a very ample divisor gives an almost Ulrich bundle.

\begin{remark}
	The existence results for Ulrich, weakly Ulrich, and almost Ulrich bundles in Theorems \ref{thm1}, \ref{thm2}, \ref{thm3}, and \ref{thm4} extend from algebraically closed fields to arbitrary fields of the stated characteristic; see Section \ref{extensionarbifields}.
\end{remark}

Combining Theorems \ref{thm1} and \ref{thm2} with the results of
\cite[Theorems 0.2 and 0.3, \S\S 2--3]{eisenbud2003resultants}
and \cite[Theorem 4.2]{eisenbud2011boij}, we obtain the following
applications to Chow forms and cohomology tables.

\begin{cor}\label{intcor}
	Let $X\subseteq\bP^N$ be a smooth projective variety over a field of
	characteristic $p>0$. A positive power of the Chow form of $X$ admits a determinantal representation.
	
	If $\dim X=2$, a positive power of the Chow form of $X$ is the determinant of a square matrix of linear forms. Moreover, the cone of
	cohomology tables of coherent sheaves on $(X,\cO_X(1))$ coincides with that on
	$(\bP^2,\cO_{\bP^2}(1))$.
\end{cor}

\section{Existence of Ulrich vector bundles on polarized smooth projective surfaces}
Let $S$ be a smooth projective surface over an algebraically closed field $k$ of characteristic $p>0$. Let $H$ be a very ample divisor on $S$. 
Throughout, stability means slope stability with respect to $H$. In this section, we show the existence of Ulrich bundles on every polarized smooth projective surface $(S, H)$

\begin{definition}
	Let $E$ be  a vector bundle on $S$.  Then $E$ is called Ulrich \wrt $H$ if the cohomology $H^\bullet(S, E(-jH)) =0$ for $1 \le j \le 2$.
\end{definition}
For other equivalent characterizations of Ulrich bundles refer to \cite[Proposition 2.1]{eisenbud2003resultants}, \cite[Theorem 2.3]{beauville2018introduction}.

The following proposition is the standard existence result obtained from the Serre construction using the Cayley-Bacharach property; see \cite[Theorem 5.1.1, Theorem 5.1.3]{huybrechts2010geometry}. We briefly recall the key steps, and the proof works in arbitrary characteristic. The construction is used in positive characteristic; see for example \cite[Proposition 7.1]{balaji2011analogue}, and also \cite[Proof of Theorem 9.5, p. 409]{langer2006moduli}.

\begin{proposition} \label{Serrecons}
	Let $L$ be a line bundle on $S$. There exists an integer $c_0$ such that, for every integer $c\ge c_0$, there is a $\mu_H$-stable rank-two vector bundle $E$ on $S$ with
	\[
	\det E\cong L,\qquad c_2(E)=c.
	\]
	In particular, the discriminant $\Delta(E)=4c_2(E)-c_1(L)^2$ can be made arbitrarily large.
\end{proposition}

\begin{proof}
	Note that it is enough to construct a stable rank $2$ bundle $G$ with $\det G=L(2aH)$ for some integer $a$ whose second Chern class $c_2(G)$ can be any sufficiently large prescribed integer. We can then take $E:=G(-aH)$.
	
	Denote  $Q:=L(2aH)$, and $M:=Q\otimes\omega_S$. Make the following choices: 
	\begin{itemize}
		\item [(1)] Choose $a$ sufficiently large that both $c_1(Q)\cdot H$ and $c_1(M)\cdot H$ are positive. 
		\item [(2)] Choose  a positive integer $b$ sufficiently large so that $bH^2>\beta$, where $\beta:=\max\big\{
		\frac{c_1(Q)\cdot H}{2},
		\ c_1(M)\cdot H
		\big\}$.
		\item [(3)] Choose a smooth integral curve $C\in |bH|$ by Bertini's theorem \cite[\href{https://stacks.math.columbia.edu/tag/0G4C}{Tag 0G4C}]{stacks-project}.
		\item [(4)] Choose an integer $\ell$ with $\ell-1>b\beta$, and let $Z=\{z_1,\ldots,z_\ell\}\subset C$ be any set of $\ell$ distinct points.
		\end{itemize}

	\textbf{Claim} No effective divisor $D$ with $D\cdot H\le\beta$ contains $Z\setminus\{z\}$ for any $z\in Z$.
	
Indeed, such a divisor cannot contain $C$ as a component, since that would imply
	\[
	D\cdot H\ge C\cdot H=bH^2>\beta.
	\]
	Thus $D$ and $C$ intersect properly, and
	\[
	D\cdot C=b(D\cdot H)\le b\beta.
	\]
	It cannot therefore contain $\ell-1>b\beta$ distinct points of $C$, proving the claim.
	
	Since $c_1(M)\cdot H\le\beta$, by the claim:
	\[
	H^0\!\left(S,\mathcal I_{Z\setminus\{z\}}\otimes M\right)=0
	\qquad\text{for every }z\in Z.
	\]
In particular, $Z$ satisfies the Cayley-Bacharach condition with respect to $M=Q\otimes\omega_S$.
	
	The Serre correspondence (see \cite[Theorem 5.1.1]{huybrechts2010geometry}) therefore produces a rank $2$ bundle: 
	\[
	0\longrightarrow\mathcal O_S
	\longrightarrow G
	\longrightarrow\mathcal I_Z\otimes Q
	\longrightarrow0,
	\]
	with  $\det G=Q$, and $c_2(G)=\ell$.

For the proof slope stability of $G$ see \cite[Theorem 5.1.3]{huybrechts2010geometry}. This argument is valid in every characteristic, including $p=2$.

Set \(E=G(-aH)\). Then $E$ is $\mu_H$-stable, and
	\[
	\det E=L,
	\qquad
	c_2(E)=\ell-a\,c_1(L)\cdot H-a^2H^2.
	\]
Since $a$, $b$, and $C$ are fixed, and every sufficiently large integer $\ell$ is allowed, $c_2(E)$ takes every sufficiently large integer value. Moreover, $\Delta(E)=\Delta(G)=4\ell-c_1(Q)^2$, which can be made arbitrarily large.
\end{proof}

For $\varepsilon \in \{0, 1\}$, set $L_{\varepsilon}:=\omega_S(\varepsilon H)$. By Proposition \ref{Serrecons}, we may fix a $\mu_H$-stable rank two vector bundle $V_\varepsilon$ with $\det V_\varepsilon\cong L_\varepsilon$. The asymptotic construction in Proposition \ref{keyprop} requires both values of $\varepsilon$ to produce Ulrich bundles \wrt $mH$ for every sufficiently large integer $m$.

Fix $\varepsilon$, and let $t$ be a sufficiently large positive integer. Set
\[
m=2t+\varepsilon,\qquad V=V_\varepsilon(tH).
\]
Then Then $V$ is $\mu_H$-stable, and $\det V=L_\varepsilon(2tH)=\omega_S(mH)$. Since  $t\gg0$,  Serre vanishing gives $H^1(S,V)=H^2(S,V)=0$. 

Put $N=h^0(S,V)$. Choose distinct points $x_1,\dots,x_N$ in $S$ and nonzero linear quotients
\[
V \twoheadrightarrow V|_{x_\nu}\twoheadrightarrow k
\]
such that the induced map
\[
\operatorname{ev}:H^0(S,V)\longrightarrow k^N
\]
is an isomorphism.  To construct these quotients, impose the conditions successively. Whenever the subspace of sections satisfying the previously imposed conditions is nonzero, choose a nonzero section in that subspace, an unused point where that section does not vanish, and a fibre quotient detecting its value. Each condition reduces the remaining dimension by exactly one. After $N$ conditions, the evaluation map is injective between two $N$-dimensional vector spaces, hence is an isomorphism.
Define
\begin{equation} \label{e201}
0\longrightarrow T\longrightarrow V \longrightarrow Q:=\bigoplus_{\nu=1}^{N}k(x_\nu) \longrightarrow0.
\end{equation}

\begin{proposition}
	The sheaf $T$ is a $\mu_H$-stable torsionfree sheaf of rank two, and
	\[
	H^i(S,T)=0\quad(i=0,1,2),\qquad
	c_1(T)=c_1(V),\qquad c_2(T)=c_2(V)+N.
	\]
	Moreover, $\det T\cong\det V\cong\omega_S(mH)$.
\end{proposition}

\begin{proof}
	Since $T\subset V$ and $\mathcal Q$ has zero-dimensional support, $T$ is torsion-free of rank two and $\det T\cong\det V$. Also, $H^i(S,\mathcal Q)=0$ for $i>0$. The long exact cohomology sequence associated with \eqref{e201} gives
	\[
	0\longrightarrow H^0(S,T) \longrightarrow H^0(S,V) \xrightarrow[\sim]{\operatorname{ev}}k^N \longrightarrow H^1(S,T) \longrightarrow0
	\]
	and $H^2(S,T)\cong H^2(S,V)=0$. Thus $H^i(S,T)=0$ for $i=0,1,2$.
	
	Every rank-one subsheaf $A\subset T$ is also a subsheaf of the stable bundle $V$, so
	\[
	\mu_H(A)<\mu_H(V)=\mu_H(T).
	\]
	Hence $T$ is $\mu_H$-stable. Finally, additivity of the Chern character in \eqref{e201} gives
	\[
	c_1(T)=c_1(V),\qquad c_2(T)=c_2(V)+N.
	\]
\end{proof}

Set $R:=T(-tH)$. Then $R$ is $\mu_H$-stable, $\det R\cong L_\varepsilon$, and
\[
\Delta(R)=\Delta(T)=\Delta(V)+4N=\Delta(V_\varepsilon)+4N.
\]
is arbitrarily large. More explicitly, since $\chi(S,T)=0$ and $c_1(T)=K_S+mH$, Riemann-Roch gives
\[
c_2(T)=2\chi(\mathcal O_S)+\frac12(K_S+mH)\cdot mH.
\]
Consequently, 
\begin{equation} \label{e202}
	\Delta(R)=8\chi(\mathcal O_S)-K_S^2+m^2H^2
\end{equation}
tends to infinity as $t$ tends to infinity, while $\det R\cong L_\varepsilon$ remains fixed.

Let $M_{L_\varepsilon}(\Delta)$ denote the moduli space of $\mu_H$-stable sheaves on $S$ with determinant $L_\varepsilon$ and discriminant $\Delta$. The moduli space is constructed as a GIT quotient $\cR//GL_{N_0}$ of a locally closed subscheme $\cR$ of the $\mathcal{Q}uot$-scheme.
The universal quotient sheaf $\mathcal U$ on $S\times\mathcal R$ is flat over $\mathcal R$.

 There exists an integer $\eta_\varepsilon$ such that for all discriminant $\Delta \ge \eta_\varepsilon$, $M_{L_\varepsilon}(\Delta)$ is an irreducible variety by \cite[Theorem 9.5]{langer2006moduli}. Moreover, after increasing $\eta_\varepsilon$ if necessary, for $\Delta \ge \eta_\varepsilon$, $M_{L_\varepsilon}(\Delta)$ contains a locally free $\mu_H$-stable sheaf. By \eqref{e202}, we may increase the lower bound on $t$ so that $\Delta(R)\ge\eta_\varepsilon$ for every $t$ under consideration. Set $\Delta:=\Delta(R)$.

\begin{lemma}
	There exists a $\mu_H$-stable rank-two vector bundle $F$ such that
	\[
	\det F\cong L_\varepsilon,\qquad
	\Delta(F)=\Delta(R),\qquad
	H^i(S,F(tH))=0\quad(i=0,1,2).
	\]
\end{lemma}

\begin{proof}
	Consider the following subsets:
	\[
	C_{\mathrm{lf}}= \{[E]\in M_{L_\varepsilon}(\Delta):E\text{ is locally free}\},
	\]
	and
	\[
	C_{\mathrm{ac},t} = \{[E]\in M_{L_\varepsilon}(\Delta):H^i(S,E(tH))=0\text{ for }i=0,1,2\}.
	\]
	Both subsets are open. Indeed, local freeness of the fibres is open in a flat family of coherent sheaves, and the required cohomology vanishings define an open subset by semicontinuity; see \cite[Theorem III.12.8]{hartshorne1977algebraic}. We may use the universal flat family $\mathcal U$ on $S\times\mathcal R$ to derive these conclusions.
	
	The set $C_{\mathrm{lf}}$ is nonempty by the choice of $\eta_\varepsilon$, and using the above mentioned result \cite[Theorem 9.5]{langer2006moduli}. The set $C_{\mathrm{ac},t}$ is also nonempty, since $R(tH)=T$ and $H^i(S,T)=0$ for $i=0,1,2$. Since $M_{L_\varepsilon}(\Delta)$ is irreducible,
	\[
	C_{\mathrm{lf}}\cap C_{\mathrm{ac},t}\ne\varnothing.
	\]
	Choose a closed point $[F]$ in this intersection. Since $k$ is algebraically closed, it is represented by a sheaf $F$ on $S$, with the required properties. 
\end{proof}
	
\begin{cor}\label{cor2}
	The bundle $W:=F(tH)$ is a $\mu_H$-stable rank-two vector bundle satisfying
	\[
	\det W\cong\omega_S(mH),\qquad
	H^i(S,W)=0\quad(i=0,1,2).
	\]
\end{cor}

\begin{proposition} \label{keyprop}
	The rank two bundle $U_m:=W(mH)$ is Ulrich with respect to $mH$.
\end{proposition}

\begin{proof}
	Since $W$ has rank two and $\det W\cong\omega_S(mH)$, we have
	\[
	W^\vee\otimes\omega_S \cong W\otimes(\det W)^{-1}\otimes\omega_S \cong W(-mH).
	\]
	By Serre duality and Corollary \ref{cor2},
	\[
	H^i(S,W(-mH)) =H^i(S, W^\vee\otimes\omega_S) \cong H^{2-i}(S,W)^\vee =0\qquad(i=0,1,2).
	\]
	Thus, for $i=0,1,2$,
	\[
	H^i(S,U_m(-mH))=H^i(S,W)=0
	\]
	and
	\[
	H^i(S,U_m(-2mH))=H^i(S,W(-mH))=0.
	\]
	Therefore $U_m$ is Ulrich with respect to $mH$.
\end{proof}

\begin{cor} \label{cor23}
	For every sufficiently large integer $m$, the polarized surface $(S,mH)$ admits a rank-two Ulrich bundle.
\end{cor}

\begin{proof}
	For each $\varepsilon\in\{0,1\}$, the construction applies to every sufficiently large integer $t$ and gives a rank two Ulrich bundle for $m=2t+\varepsilon$. Taking a common lower bound for both the even and odd parities proves the assertion.
\end{proof}

\begin{proof} [Proof of Theorem \ref{thm1}]
	Let $F: S \longrightarrow S$ be the absolute Frobenius morphism. For $q=p^e$ sufficiently large, choose a rank two Ulrich bundle $U_q$ with respect to $qH$ by Corollary \ref{cor23}, and set
	\[
	E_e:=(F^e)_*U_q.
	\]
	Since $S$ is smooth, the morphism $F^e$ is finite flat of degree $p^{2e}=q^2$. Thus $E_e$ is locally free of rank $2q^2$.
	
	For $j=1,2$, the projection formula gives
	\[
	E_e(-jH) \cong (F_S^e)_*\!\left(U_q\otimes(F_S^e)^*\mathcal O_S(-jH)\right) \cong (F_S^e)_*\bigl(U_q(-jqH)\bigr).
	\]
	Since $F^e$ is finite, it follows that
	\[
	H^i(S,E_e(-jH))\cong H^i(S,U_q(-jqH))=0
	\]
	for $i=0,1,2$ and $j=1,2$. 
	Hence $E_e$ is Ulrich with respect to $H$.
\end{proof}

\begin{remark} \label{specialul}
	We have $\det U_q \simeq \det W \otimes \cO_S(2qH) \simeq \omega_S(3qH)$. For $q=p^e$, $E_e=F^e_*U_q$ is Ulrich \wrt $H$ by the preceding proof.  Since $U_q$ has rank $2$, the isomorphism $(F^e)^!\omega_S\simeq\omega_S$, finite duality \cite[Exercise III.6.10(b), p. 239]{hartshorne1977algebraic}, and the projection formula give
	\begin{equation} \label{fduality}
		E_e^\vee\otimes\omega_S \simeq \cH om(F^e_* U_q, \omega_S) \simeq F^e_*(U_q^\vee\otimes\omega_S)
		\simeq F^e_*\bigl(U_q(-3qH)\bigr) \simeq E_e(-3H).
	\end{equation}
	Thus, $E_e\simeq E_e^\vee\otimes\omega_S(3H)$.  Moreover, the composite
	\[
	\beta:E_e\otimes E_e \longrightarrow F^e_*(\wedge^2U_q) \simeq F^e_*(\omega_S \otimes (F^e)^* \cO_S(3H)) \simeq F^e_*\omega_S\otimes \cO_S(3H)  \xrightarrow{\operatorname{Tr}_{F^e}\otimes1}
		\omega_S\otimes \cO_S(3H),
	\]
	where $\operatorname{Tr}_{F^e}$ denotes the Grothendieck trace, is a perfect alternating pairing. Indeed, $\beta(u\otimes u)=0$ because $u\wedge u=0$, so $\beta$ is alternating. The pairing is perfect by \eqref{fduality}.
	
Thus, $E_e$ is skew symmetric in the sense of \cite[p. 553]{eisenbud2003resultants}. Furthermore, since $r:=\rk(E_e)=2q^2$, the perfect alternating pairing gives $\det E_e\simeq \bigl(\omega_S(3H)\bigr)^{\otimes (r/2)}$.
Thus, $E_e$ is a special Ulrich bundle \wrt $H$ \cite[p. 2]{casnati2020examples}.

Since $\rk(E_e)$ is even, for $p>2$, it follows that some positive power of Chow forms is Pfaffian of a skew symmetric matrix with entries linear in the Pl\"{u}cker coordinates; see \cite[Theorem 0.3, Theorem 3.3]{eisenbud2003resultants}. 
	This sharpens the linear determinantal representation of Chow forms in Corollary \ref{intcor}. 
	\end{remark}

\begin{remark}
	Set $r:=\rk (E_e)=2q^2$. Since $E_e$ is Ulrich, $\chi(E_e)=rH^2$. By Remark \ref{specialul}, we have
	\[
	c_1(E_e)=\frac r2(K_S+3H).
	\]
	By the Riemann-Roch theorem, we obtain
	\[
	c_2(E_e)=r\chi(\mathcal O_S)+\frac12c_1(E_e)^2-\frac12c_1(E_e)K_S-rH^2.
	\]
	This gives
	\[
	\frac{\Delta(E_e)}{r^2}=\frac{2r\,c_2(E_e)-(r-1)c_1(E_e)^2}{r^2}=\frac{H^2+8\chi(\mathcal O_S)-K_S^2}{4}.
	\]
	Therefore, whenever $H^2<K_S^2-8\chi(\mathcal O_S)$, $\Delta(E_e)<0$, and thus $E_e$ cannot be strongly semistable \wrt $H$; see \cite[Corollary 3.12]{langer2022boundedness}.
\end{remark}

\subsection{Examples arising from the surface case}
\label{productcons}
Let $X_1,X_2$ be smooth projective varieties over an algebraically closed field of arbitrary characteristic. Let $d_i=\dim X_i\geq1$, for $i=1, 2$.
Let $H_i$ be a very ample divisor on $X_i$. Suppose that $E_i$ is an Ulrich bundle on $(X_i,H_i)$.
Put
\[
X=X_1\times X_2,\qquad \mathcal O_X(H)= \mathcal O_{X_1}(H_1)\boxtimes \mathcal O_{X_s}(H_2).
\]
Then the vector bundle $\mathcal U=E_1\boxtimes E_2(d_1H_2)$ is Ulrich with respect to $H$; see \cite[Proposition 2.6(b)]{eisenbud2003resultants}, and \cite[property (2.5), p. 4]{beauville2018introduction}. 

By induction, this construction extends to finite products. Let $X_1,\ldots,X_s$ be smooth projective varieties of dimensions $d_1, \cdots, d_s$. Let $H_i$ be a very ample divisor on $X_i$. Suppose that $E_i$ is an Ulrich bundle on $(X_i,H_i)$. Put
\[
X=X_1\times\cdots\times X_s,\qquad \mathcal O_X(H)= \mathcal O_{X_1}(H_1)\boxtimes\cdots\boxtimes \mathcal O_{X_s}(H_s).
\]
Then the vector bundle $\mathcal U=E_1 \boxtimes E_2(d_1H_2) \boxtimes E_3((d_1+d_2)H_3)
\boxtimes \cdots\boxtimes E_s((d_1+d_2+\cdots+d_{s-1})H_s)$ is Ulrich on $(X, H)$.
Combining this product construction with the existence Theorem \ref{thm1} gives the following consequence.

\begin{cor}
	Let $k$ be an algebraically closed field of characteristic
	$p>0$. Each of the following polarized varieties admits an
	Ulrich bundle:
	\begin{enumerate}
		\item $X=S_1\times\cdots\times S_m$, where $S_i$ is a smooth projective surface
		over $k$, with polarization $	\mathcal O_X(H)=
		\mathcal O_{S_1}(H_1)\boxtimes\cdots\boxtimes
		\mathcal O_{S_m}(H_m)$ for very ample divisor $H_i$.
		
		\item $X=S\times\mathbb P^n$, where $S$ is a smooth
		projective surface over $k$, with polarization $\mathcal O_X(H)=
		\mathcal O_S(H_S)\boxtimes
		\mathcal O_{\mathbb P^n}(a)$, where $H_S$ is very ample and $a\geq 1$ is an integer.
		
		\item $X=S\times C$, where $S$ is a smooth projective
		surface and $C$ is a smooth projective curve over $k$, with
		polarization $\mathcal O_X(H)=
		\mathcal O_S(H_S)\boxtimes\mathcal O_C(H_C)$, where $H_S$ and $H_C$ are very ample.
	\end{enumerate}
\end{cor}

\begin{proof}
	Assertion (1) follows from combining Theorem \ref{thm1} and the above discussion.
	
	For (2), let $E$ be an Ulrich bundle on $(S,H_S)$, and let
	$V_a$ be an Ulrich bundle on
	$(\mathbb P^n,\mathcal O_{\mathbb P^n}(a))$; see Proposition \ref{ulproj}. Then $E\boxtimes V_a(2a)$
	is an Ulrich bundle on $X$.
	
	For (3) Let $g$ be the genus of $C$. Any line bundle of the form $M=F \otimes O_C(H_C)$ is Ulrich, where $F \in \text{Pic}^{g-1}(C) \setminus \Theta$, and $\Theta$ is the theta divisor.
	Choosing Ulrich $E$ on $(S,H_S)$, $E \boxtimes F(3H_C)$ gives an Ulrich bundle on $(X, H)$.
\end{proof}

We next apply Theorem \ref{thm1} to projective bundles over surfaces.
We first extend the required cohomological vanishings to ample divisors.

\begin{lemma} \label{lemample}
	Let $S$ be a smooth projective surface over an algebraically closed
	field $k$ of characteristic $p>0$, and let $A$ be an ample divisor on
	$S$. For every sufficiently large integer $e$, there exists a vector
	bundle $C_e$ of rank $2p^{2e}$ such that
	\[
	H^i(S,C_e(-jA))=0
	\qquad (i=0,1,2,\quad j=1,2).
	\]
\end{lemma}

\begin{proof}
	Choose an integer $b\geq0$ such that $p^bA$ is very ample. By
	Theorem \ref{thm1}, for every sufficiently large integer $t$, there
	exists a vector bundle $G_t$ of rank $2p^{2t}$ that is Ulrich with
	respect to $p^bA$. Let $F_S:S\longrightarrow S$ be the absolute
	Frobenius morphism, and set
	\[
	C_{t+b}:=(F_S^b)_*G_t.
	\]
	Since $F_S^b$ is finite flat of degree $p^{2b}$, the sheaf $C_{t+b}$
	is locally free of rank $2p^{2(t+b)}$. For $i=0,1,2$ and $j=1,2$,
	the projection formula gives
	\[
	H^i(S,C_{t+b}(-jA))
	\cong H^i(S,G_t(-jp^bA))=0.
	\]
	Relabelling $e=t+b$ proves the assertion.
\end{proof}

Let $V$ be a vector bundle of rank $r \ge 2$ on $S$. We consider the projective bundle
\[
\pi:X=\mathbb P_S(V)=\operatorname{Proj}_S(\operatorname{Sym}V) \longrightarrow S.
\]
Thus $\mathbb P_S(V)$ parametrizes one dimensional quotients of $V$, and there is a tautological surjection
\[
\pi^*V\twoheadrightarrow \mathcal O_X(1),
\]
where $\mathcal O_X(1)$ is the tautological line bundle. With this convention, we have $\pi_*\mathcal O_X(1)\cong V$. Choose a divisor $D$ on $S$ such that
$\cO_S(D)\cong\det V$. For a vector bundle $E$ on $S$, set $\mathcal U:=\pi^*E\otimes\cO_X(1)$.

\begin{proposition} 
	\label{hochenegger}
	Let $S$ be a smooth projective surface over an algebraically closed
	field of arbitrary characteristic. 
	Assume that the tautological line bundle $\cO_X(1)$ is very ample. 
Then $H^\bullet(X,\mathcal U(-j))=0$ for $1\leq j\leq r+1$  if and only if
\[
H^\bullet(S,E)=H^\bullet(S,E(-D))=0.
\]
\end{proposition}

\begin{proof}
	Since $\dim X=r+1$, this is
	\cite[Theorem 5.1]{hochenegger2022ulrich}, with $A=0$ in the
	notation of that paper. Although the paper assumes characteristic
	zero, its proof via \cite[Proposition 3.1]{hochenegger2022ulrich}
	remains valid in arbitrary characteristic: it uses the standard
	semiorthogonal decomposition for projective bundles, the canonical
	bundle formula, Serre duality, and the derived adjunction
	$L\pi^*\dashv R\pi_*$, all of which hold in arbitrary characteristic.
	\end{proof}

\begin{theorem} \label{projbunthm}
	Let $S$ be a smooth projective surface over an algebraically closed
	field $k$ of characteristic $p>0$. 
	Let $\pi:X=\bP_S(V)\longrightarrow S$ be the projective bundle associated to a vector bundle  $V$
of rank $r \geq2$ on $S$. Let $\cO_X(1)$ be the tautological line bundle on $X$. 
	Let $H$ be a very ample divisor on $X$ such that $\cO_X(H)\cong\cO_X(a)\otimes\pi^*M$,
	where $a\geq1$ is an integer and $M\in\operatorname{Pic}(S)$.
	\begin{enumerate}
		\item[(1)] If $a=1$, then $(X,H)$ admits an Ulrich bundle of rank
		$2p^{2e}$ for every sufficiently large integer $e$.
		\item[(2)] If $V$ is a direct sum of line bundles, then there exists
		an integer $\delta\geq1$, independent of $e$, such that $(X,H)$
		admits an Ulrich bundle of rank $2\delta(r+1)!p^{2e}$ for every
		sufficiently large integer $e$.
	\end{enumerate}
\end{theorem}

\begin{proof}
	\textup{(1)} Set $W:=V\otimes M$. Under the natural identification
	$\bP_S(W)\cong X$, the line bundle $\cO_X(H)$ is tautological.
	Hence $W$ is ample, and so is $\det W$ by
	\cite[Proposition 3.2 and Corollary 2.6]{hartshorne1966ample}.
	Choose an ample divisor $D$ on $S$ such that
	$\cO_S(D)\cong\det W$.
	
	By Lemma \ref{lemample}, for every sufficiently large integer $e$,
	there exists a vector bundle $C_e$ of rank $2p^{2e}$ satisfying
	\[
	H^\bullet(S,C_e(-D))=H^\bullet(S,C_e(-2D))=0.
	\]
Applying Proposition \ref{hochenegger} to $\bP_S(W)$ and $C_e(-D)$ shows that
\begin{equation}
		E_e:=\pi^*(C_e(-D))\otimes\cO_X(H)
	\end{equation}
	is an Ulrich bundle on $(X,H)$ of rank $2p^{2e}$.
	
	\smallskip
	\noindent\textup{(2)} Write $V=\bigoplus_{\nu=1}^rL_\nu$.
	By  \cite[Lemma 4.5]{langer2022boundedness}, there exist a smooth
	projective surface $T$, a finite surjective morphism
	$f:T\longrightarrow S$, and a line bundle $N$ on $T$ such that
	\[
	N^{\otimes a}\cong f^*M.
	\]
	Since $T$ and $S$
	are smooth surfaces, $f$ is flat; see \cite[\tagref{00R4}]{stacks-project}.
	Set
	\[
	W:=f^*V\otimes N =\bigoplus_{\nu=1}^r A_\nu,
	\qquad
	A_\nu:=f^*L_\nu\otimes N.
	\]
	Under the isomorphism $\bP_T(W)\cong X\times_S T$, the pullback of
	$\cO_X(H)$ is $\cO_{\bP_T(W)}(a)$. Thus $W$ is ample,
	and hence, each quotient line bundle $A_\nu$ is ample.
	
	Let $F_T:T\longrightarrow T$ be absolute Frobenius.
	Choose $b\geq0$ such that every $A_\nu^{\otimes p^b}$ is
	very ample, and put
	\[
	f_b:=f\circ F_T^b,\qquad N_b:=(F_T^b)^*N,\qquad W_b:=(F_T^b)^*W
	\cong f_b^*V\otimes N_b.
	\]
	Consider the projective bundle $\rho:Y=\bP_T(W_b)\longrightarrow T$.
	Since each $A_\nu^{\otimes p^b}$ is very ample, the vector bundle
	$W_b\cong\bigoplus_{\nu=1}^r A_\nu^{\otimes p^b}$
	is very ample by \cite[Proposition 1.2(2)]{sumihiro1982theorem}.
	Hence its tautological line bundle
	$\cO_Y(\eta)=\cO_{\bP_T(W_b)}(1)$ is very ample.
Base change along $f_b$ gives a finite flat
	morphism $g:Y\longrightarrow X$ such that
	\begin{equation}\label{pbe1}
		g^*\cO_X(H)\cong\cO_Y(a\eta),
	\end{equation}
	because $N_b^{\otimes a}\cong f_b^*M$.
	Fix these choices and set $\delta:=\deg(g)=(\deg f)p^{2b}$.
	
	Put $n:=\dim Y=r+1$. By \textup{(1)}, for every sufficiently
	large integer $e$, there exists an Ulrich bundle $\mathcal U_e$
	on $(Y,\eta)$ of rank $2p^{2e}$. Choose a finite linear projection
	$\lambda:Y\longrightarrow\bP^n$ associated with $\eta$.
	By Proposition \ref{ulproj}, there exists an Ulrich
	bundle $V_a$ of rank $n!$ on $(\bP^n,\cO_{\bP^n}(a))$. Then
	\[
	\widetilde{\mathcal U}_e:=\mathcal U_e\otimes\lambda^*V_a
	\]
	is Ulrich with respect to $a\eta$ by \cite[Proposition 5.4]{eisenbud2003resultants}.
	
	Set $E_e:=g_*\widetilde{\mathcal U}_e$.
	Since $g$ is finite flat, $E_e$ is locally free of rank $\operatorname{rk}(E_e) =\delta\,\operatorname{rk}(\widetilde{\mathcal U}_e)=2\delta(r+1)!p^{2e}$. By the projection formula and \eqref{pbe1},
	\[
	H^i(X,E_e(-jH)) \cong H^i(Y,\widetilde{\mathcal U}_e(-ja\eta)) =0
	\qquad
	(0\leq i\leq n,\quad 1\leq j\leq n).
	\]
	Thus $E_e$ is Ulrich with respect to $H$, and $\delta$
	is independent of $e$.
\end{proof}

\begin{remark}
The splitting assumption in Theorem \ref{projbunthm}\textup{(2)}
can be replaced by the existence of a filtration by subbundles
			\[
			0=V_0\subset V_1\subset\cdots\subset V_r=V
			\]
			with line bundle quotients $L_\nu:=V_\nu/V_{\nu-1}$ such that
			$L_\nu^{\otimes a}\otimes M$ is ample for every $\nu$.
			
			Indeed, after the finite cover $f:T\longrightarrow S$ used in the proof, $W=f^*V\otimes N$ has a filtration with ample line bundle
			quotients $A_\nu=f^*L_\nu\otimes N$.
			Fix a very ample line bundle $B$ on $T$. For sufficiently large $b$,
			each $A_\nu^{\otimes p^b}\otimes B^{-1}$ is globally generated and has vanishing first cohomology. Induction along the 
			filtration therefore shows that $Q_b:=(F_T^b)^*W\otimes B^{-1}$ is globally generated.
			 Then it follows that $(F_T^b)^*W=Q_b\otimes B$ is very ample.
			The remainder of the proof applies.
		\end{remark}

\section{Weakly Ulrich bundles on polarized smooth projective varieties} \label{weaklysec}
Let $X$ be a smooth projective variety over an algebraically closed field $k$, and let $H$ be a very ample divisor on $X$. In this section, we construct weakly Ulrich bundles on $X$ with respect to $qH$ for all sufficiently large integers $q$. As a consequence, we prove that $(X,H)$ admits weakly Ulrich bundles when $k$ has positive characteristic.

We recall the notion of weakly Ulrich bundles from \cite[Section 2]{eisenbud2003resultants}; see also \cite[Section 2]{Jos21}.
For a coherent sheaf $E$ and an integer $t$, write $E(tH)=E\otimes\cO_X(tH)$. 

\begin{definition}\label{defweak}
	A vector bundle $E$ on a polarized smooth projective $n$-fold $(X,H)$ is \emph{weakly Ulrich} if
	\begin{align}
		H^i(X,E(-mH))&=0 &&(1\leq i\leq n,\ m\leq i-1),\label{weak1}\\
		H^i(X,E(-mH))&=0 &&(0\leq i\leq n-1,\ m\geq i+2).\label{weak2},
	\end{align}
	for all integers $m$ in the indicated ranges. 
\end{definition}

For our construction, it is convenient to fix the twist parameter
$m$ and list the required cohomological vanishings. We use the
following equivalent formulation in Proposition \ref{prop33}.

\begin{lemma}\label{equiweaklem}
	Let $E$ be a vector bundle on $(X,H)$. Then $E$ is weakly Ulrich
	if and only if $H^i(X,E(-mH))=0$ for every integer $m$ and every
	degree $i$ in the ranges specified below:
	\[
	\begin{array}{c|l}
		\text{Range of }m
		& \text{Required vanishing range for }i \\ \hline
		m\leq 0
		& 1\leq i\leq n \\[2pt]
		1\leq m\leq n
		& 0\leq i\leq m-2
		\quad\text{or}\quad m+1\leq i\leq n \\[2pt]
		m\geq n+1
		& 0\leq i\leq n-1.
	\end{array}
	\]
	An empty range imposes no condition.
\end{lemma}

\begin{proof}
	For $1\leq m\leq n$, condition \eqref{weak1}
	requires vanishing in degrees $i\geq m+1$, while
	condition \eqref{weak2} requires vanishing in degrees
	$i\leq m-2$. For $m\leq 0$, the first condition requires
	vanishing in every positive degree. For $m\geq n+1$, the
	second condition requires vanishing in every degree below $n$.
	Thus the table lists exactly the vanishings in
	Definition \ref{defweak}, grouped according to $m$.
\end{proof}

\subsection{Cohomology on a product of projective lines} \label{prodcoho}

Throughout this section, $k$ is assumed to be an algebraically closed field of arbitrary characteristic. For $n \ge 1$, we consider
\[
Z=(\bP^1_k)^n,\qquad
\cO_Z(\avec)=\boxtimes_{\nu=1}^n\cO_{\bP^1}(a_\nu),\qquad
\one=(1,\ldots,1),
\]
where $\avec$ is the $n$-tuple $(a_1, a_2, \cdots, a_n)$. For a sheaf $G$ on $Z$, denote $G(\avec)=G\otimes\cO_Z(\avec)$.

\begin{lemma} \label{lem31}
	Let $\{1,\ldots,n\}=P\sqcup N\sqcup J$ be a partition, and put $r=|N|$ and $s=|J|$. Suppose that
	\[
	a_\nu\geq0\ (\nu\in P),\qquad a_\nu\leq-2\  (\nu\in N),
	\]
	with no restriction on $a_{\nu}$ for $\nu\in J$. Then $H^i(Z,\cO_Z(\avec))=0$ unless $r\leq i\leq r+s$.
\end{lemma}
\begin{proof}
	By the K\"unneth formula, each summand of
	$H^i(Z,\cO_Z(\avec))$ has the form
	\[
	\bigotimes_{\nu=1}^n H^{i_\nu}(\bP^1,\cO_{\bP^1}(a_\nu)),
	\qquad i_\nu\in\{0,1\},\quad \sum_{\nu=1}^n i_\nu=i;
	\]
	see \cite[\tagref{0BED}]{stacks-project}. For such a summand to be nonzero,
	one must have $i_\nu=0$ for $\nu\in P$ and $i_\nu=1$ for
	$\nu\in N$, by the cohomology of line bundles on $\bP^1$.
	The coordinates in $J$ contribute between $0$ and $s$ to the total
	degree. Thus $r\leq i\leq r+s$.
\end{proof}

\begin{lemma}\label{lem32}
	For every vector bundle $E$ on $Z$, there are vector bundles $K_0,\ldots,K_n$, with $K_0=E$, and exact sequences
	\begin{equation}\label{eq312}
		0\longrightarrow K_{j+1}\longrightarrow P_j \longrightarrow K_j \longrightarrow0 \qquad(0\leq j <n),
	\end{equation}
	where $P_j=\cO_Z(-b_j \one)^{\oplus t_j}$ for some integer
	$b_j > 0$ and some integer $t_j \ge 0$.
\end{lemma}
\begin{proof}
	Set $K_0=E$ and proceed inductively. Since $\cO_Z(\one)$ is ample,
	choose $b_j>0$ such that $K_j(b_j\one)$ is globally generated.
	The evaluation map gives a surjection
	\[
	\cO_Z^{\oplus t_j}\twoheadrightarrow K_j(b_j\one)
	\]
	for some integer $t_j\geq0$. Twisting by $\cO_Z(-b_j\one)$ gives
	a surjection $P_j\twoheadrightarrow K_j$. Let $K_{j+1}$ be its
	kernel. Since $K_j$ is locally free, the surjection splits locally, so $K_{j+1}$ is locally free.
\end{proof}

\begin{proposition}\label{prop31}
	Let $E$ be a vector bundle on $Z$. There is an integer $c_E\geq2$
	such that, for every partition $\{1,\ldots,n\}=P\sqcup N\sqcup J$
	and every multidegree $\avec=(a_1, \cdots, a_n) \in\bZ^n$ satisfying
	\[
	a_\nu\geq c_E\quad(\nu\in P),\qquad
    a_\nu\leq-c_E\quad(\nu\in N),
	\]
	with no restriction on $a_\nu$ for $\nu\in J$, one has
	$H^i(Z,E(\avec))=0$ for all $i>|N|+|J|$.
	\end{proposition}
\begin{proof}
	Choose the $n$ exact sequences given by Lemma \ref{lem32}, and take
	\[
	c_E\geq\max\{2,b_0+2,\ldots,b_{n-1}+2\}.
	\]
	Every summand of $P_j(\avec)$ is $\cO_Z(\avec-b_j\one)$.
	For $\nu\in P$, we have
	$a_\nu-b_j\geq c_E-b_j\geq0$, whereas for $\nu\in N$, we have
	$a_\nu-b_j\leq-c_E-b_j\leq-2$. The coordinates indexed by $J$
	are unrestricted. Hence Lemma \ref{lem31} gives
	\[
	H^i(Z,P_j(\avec))=0\qquad(i>|N|+|J|).
	\]
	Fix $i>|N|+|J|$. Since $H^{i+j}(Z,P_j(\avec))=0$ for
	$0\leq j<n$, the long exact sequences associated to
	\eqref{eq312}, after twisting by $\cO_Z(\avec)$, give
	injections
	\[
	H^i(Z,E(\avec)) \hookrightarrow H^{i+1}(Z,K_1(\avec))
	\hookrightarrow \cdots\hookrightarrow H^{i+n}(Z,K_n(\avec)).
	\]
	Since $i>0$ and $\dim Z=n$, the last group vanishes. Thus
	$H^i(Z,E(\avec))=0$.
\end{proof}

\begin{proposition} \label{prop2}
	Let $E$ be a vector bundle on $Z$. There is an integer $c \geq2$ such that, for every partition $\{1,\ldots,n\}=P\sqcup N\sqcup J$
	and every multidegree $\avec$ satisfying
	\[
	a_\nu\geq c\ (\nu\in P),\qquad
	a_\nu\leq-c\ (\nu\in N),
	\]
	with no restriction for $\nu\in J$, one has
	\begin{equation} 
		H^i(Z,E(\avec))=0
		\quad\text{unless}\quad |N| \leq i \leq |N|+|J|.
	\end{equation}
\end{proposition}
\begin{proof}
	Apply Proposition \ref{prop31} to both $E$ and $E^\vee\otimes\omega_Z$, and take $c=\max \{c_E, c_{E^\vee\otimes\omega_Z}\}$.
	 Hence, Proposition \ref{prop31}  gives $H^i(Z,E(\avec))=0$ for $i > |N|+|J|$.
	
	For the lower bound, the multidegree $-\avec$ has $|P|$ sufficiently negative coordinates and $|J|$ unrestricted coordinates. Hence Proposition \ref{prop31}  gives
	\[
	H^j(Z,E^\vee\otimes\omega_Z(-\avec))=0
	\qquad(j>|P|+|J|=n-|N|).
	\]
	If $i<|N|$, then $n-i>n-|N|$, so Serre duality gives
	\[
	H^i(Z,E(\avec))
	\cong H^{n-i}\bigl(Z,(E^\vee\otimes\omega_Z)(-\avec)\bigr)^\vee
	=0.
	\] 
\end{proof}

\subsection{Asymptotic existence of weakly Ulrich bundles} \label{asympweak}

We continue to work over an algebraically closed field $k$ of arbitrary
characteristic. Let $X$ be a smooth projective variety of dimension
$n\geq1$ over $k$, and let $H$ be a very ample divisor giving an
embedding $X \hookrightarrow \bP^N$. Since $k$ is infinite, we may
choose a linear projection with centre disjoint from the embedding of $X$ that induces
a finite surjective morphism
\[
\pi:X\longrightarrow\bP^n
\]
such that $\pi^*\cO_{\bP^n}(1)\cong\cO_X(H)$; see
\cite[Proposition 2.13]{IMW23}. Its degree is $H^n$. Since $X$ and
$\bP^n$ are smooth of the same dimension, $\pi$ is flat; see \cite[\tagref{00R4}]{stacks-project}.

Now we consider a finite map
\[
\rho:Z=(\bP^1)^n\longrightarrow \bP^n.
\]
Such a map exists by sending $(x_1, x_2, \cdots, x_n) \in (\bP^1)^n$ to its class in $Sym^n(\bP^1)\cong \bP^n$. The map $\rho$ is finite flat of degree $n!$, and $\rho^*\cO_{\bP^n}(1)=\cO_Z(\one)$; see for example \cite[Theorem 6.1]{eisenbud2009betti}. 

Fix $\pi$ and $\rho$, and set
\begin{equation}\label{eq325}
	A=\pi_*\cO_X,\qquad G=\rho^*A.
\end{equation}
Since $\pi$ is finite flat, $A$ and $G$ are vector bundles. For every positive integer $q$, define
\begin{equation}\label{eq210}
	L_q=\cO_Z(q-1,2q-1,\ldots,nq-1),\qquad
	V_q=\rho_*L_q,\qquad B_q=\pi^*V_q.
\end{equation}
Since $\rho$ is finite flat of degree $n!$, $V_q$ and $B_q$ are vector bundles
of rank $n!$.

The bundles $V_q$ give nice examples of Ulrich bundles on $\bP^n$.
\begin{proposition}\label{ulproj}
	For every integer $q>0$, the bundle $V_q$ on $\bP^n$ is Ulrich \wrt $\cO_{\bP^n}(q)$.
\end{proposition}
\begin{proof}
	The line bundle $L_q$ is Ulrich \wrt $\cO_Z(q\one)$, and thus, $V_q$ is Ulrich \wrt $\cO_{\bP^n}(q)$; for more details see \cite[Proposition 4]{beauville2018introduction}.
\end{proof}

We now express the cohomology of $B_q$ in terms of the fixed bundle $G$.

\begin{lemma}\label{lem321}
	For every integer $m$ and every $i \ge 0$, there is an isomorphism
	\begin{equation}\label{eq:cohomology-identity}
		H^i(X,B_q(-mqH)) \cong H^i\bigl(Z,G\otimes \cO_Z((1-m)q-1,\ldots,(n-m)q-1)\bigr).
		\end{equation}
\end{lemma}
\begin{proof}
	We have $B_q(-mqH) \cong \pi^*V_q \otimes \pi^* \cO_{\bP^n}(-mq) \cong \pi^* V_q(-mq)$. Recall that $V_q=\rho_*L_q$. Since $\pi$ and
	$\rho$ are finite, applying the projection formula to each gives
	\begin{align*}
		H^i(X,B_q(-mqH))
		&\cong H^i(\bP^n, \pi_* \cO_X \otimes V_q(-mq))\\
		& \cong H^i(Z, \rho_* \big(\rho^*A \otimes\rho^*\cO_{\bP^n}(-mq) \otimes L_q\big))\\
		&\cong H^i(Z,\rho^*A\otimes L_q\otimes \rho^*\cO_{\bP^n}(-mq)).
	\end{align*}
	Using $\rho^* A=G$, $\rho^*\cO_{\bP^n}(-mq) \cong \cO_Z(-mq \one)$, and the definition of $L_q$ from \eqref{eq210}, the result follows.
\end{proof}

\begin{proposition}\label{prop33}
	There is an integer $q_0 \geq1$, independent of the twist parameter
	$m$, such that, for every integer $q\geq q_0$,
	\begin{align}
		H^i(X,B_q(-mqH))&=0 &&(m\leq0,\ i>0),\label{eq:pos-range}\\
		H^i(X,B_q(-mqH))&=0 &&(1\leq m\leq n,\ i\notin\{m-1,m\}),\label{eq:middle-range}\\
		H^i(X,B_q(-mqH))&=0 &&(m\geq n+1,\ i<n).\label{eq:neg-range}
	\end{align}
\end{proposition}
\begin{proof}
	Let $c$ be the constant in Proposition \ref{prop2} for the fixed
	bundle $G$ in \eqref{eq325}, and take $q_0=c+1$.
	Fix $q\geq q_0$. By Lemma \ref{lem321}, it
	suffices to apply Proposition \ref{prop2} to the multidegree
	$\avec$ with $a_\nu=(\nu-m)q-1$. We distinguish three cases.
	
	\smallskip\noindent
	\emph{Case 1: $m\leq0$.}
	Every coordinate satisfies $a_\nu\geq q-1\geq c$.
	Take $P=\{1,\ldots,n\}$ and $N=J=\varnothing$. Thus, the interval $[|N|, |N|+|J|]$ is $[0,0]$, which gives $H^i(X,B_q(-mqH))=0$ for $i>0$.
	
	\smallskip\noindent
	\emph{Case 2: $1\leq m\leq n$.}
	Here the coordinates satisfy
	\[
	\begin{cases}
		a_\nu\leq-q-1\leq-c,&\nu<m,\\
		a_m=-1,&\nu=m,\\
		a_\nu\geq q-1\geq c,&\nu>m.
	\end{cases}
	\]
	Take $N=\{1,\ldots,m-1\}$, $J=\{m\}$, and
	$P=\{m+1,\ldots,n\}$. Since $|N|=m-1$ and $|J|=1$,
	Proposition \ref{prop2} gives $H^i(X,B_q(-mqH))=0$ unless $i\in\{m-1,m\}$.
	
\smallskip\noindent
	\emph{Case 3: $m\geq n+1$.}
	Every coordinate satisfies $a_\nu\leq-q-1\leq-c$. Taking
	$N=\{1,\ldots,n\}$ and $P=J=\varnothing$ gives the interval $[n,n]$. Thus, $H^i(X,B_q(-mqH))=0$ for all $i<n$
	\end{proof}

\begin{proof}[Proof of Theorem \ref{thm3}]
	The bundle $B_q$ defined in \eqref{eq210} has rank $n!$.
	By Proposition \ref{prop33}, for every $q\geq q_0$ its cohomology
	relative to the polarization $qH$ vanishes in the
	required ranges of Lemma \ref{equiweaklem}. Hence $B_q$ is
	weakly Ulrich with respect to $qH$.
	\end{proof}

\subsection{Existence of weakly Ulrich for fixed polarization in positive characteristic}
\label{frobdirectweak}

In this subsection, we assume that $k$ is algebraically closed field of
characteristic $p>0$. Let $X$ be a smooth projective variety of dimension $n$ over $k$. We consider absolute Frobenius morphism $F:X \longrightarrow X$. The morphism $F$ is finite of degree $p^{n}$. Since $X$ smooth, $F$ is flat.  Consequently its $e$-th iterate $F^e=F \circ F \circ \cdots \circ F:X \longrightarrow X$ is finite flat of degree $p^{en}$.

\begin{proposition}\label{prop34}
	Let $e\geq1$, put $q=p^e$, and let $B$ be a vector bundle on $X$.
	Let $H$ be a very ample divisor on $X$. Then $F^e_*B$ is locally
	free of rank $q^n\operatorname{rk}(B)$. Moreover, $B$ is weakly
	Ulrich (resp. Ulrich) with respect to $qH$ if and only if
	$F^e_*B$ is weakly Ulrich (resp. Ulrich) with respect to $H$.
\end{proposition}
\begin{proof}
	Since $F^e$ is finite flat of degree $q^n$, the sheaf $F^e_*B$
	is locally free of rank $q^n\operatorname{rk}(B)$.
	For every integer $m$, we have
	$(F^e)^*\cO_X(-mH)\cong\cO_X(-mqH)$, so the projection
	formula gives
	\[
	(F^e_*B)(-mH)
	\cong F^e_*\bigl(B\otimes(F^e)^*\cO_X(-mH)\bigr) \cong F^e_*(B(-mqH)).
	\]
	Since $F^e$ is finite, it follows that
	\begin{equation}\label{eq320}
		H^i(X,(F^e_*B)(-mH)) \cong H^i(X,B(-mqH)).
	\end{equation}
	Applying \eqref{eq320} in the respective
	vanishing ranges proves both equivalences.
\end{proof}

\begin{proof}[Proof of Theorem \ref{thm2}]
	Take $q_0$ from Proposition \ref{prop33}, and choose $e_0\geq1$
	such that $p^{e_0}\geq q_0$. For every $e\geq e_0$, the bundle
	$B_{p^e}$ is weakly Ulrich with respect to $p^eH$ by
	Theorem \ref{thm3}. By Proposition \ref{prop34}, the bundle
	\[
	E_e=F^e_*B_{p^e} =F^e_*\bigl(\pi^*\rho_* \cO_Z(p^e-1,2p^e-1,\ldots,np^e-1)\bigr)
	\]
	is weakly Ulrich with respect to $H$ and has rank $n!p^{en}$.
\end{proof}

\begin{proposition}\label{initialized}
	For every integer $q\geq1$, the bundle $B_q$ is initialized
	with respect to $qH$. Consequently, for every $e\geq1$,
	the bundle $E_e=(F_X^e)_*B_{p^e}$ is initialized with
	respect to $H$.
\end{proposition}

\begin{proof}
	Fix $q\geq1$, and set $\avec=(-1,q-1,\ldots,(n-1)q-1)$, and $M=\cO_Z(\avec)$. By Lemma \ref{lem321},
	\[
	H^0(X,B_q(-qH))\cong H^0(Z,G\otimes M).
	\]
	Let $c$ be the constant in Proposition \ref{prop2} for $G$.
	Applying that proposition to $c\avec$, with
	$N=\{1\}$, $J=\{2,\ldots,n\}$, and $P=\varnothing$, gives
	\[
	H^0(Z,G\otimes M^{\otimes c})=0.
	\]
	
	Multiplication in the $\cO_Z$-algebra $G=\rho^*\pi_*\mathcal O_X$ induces a map
	\[
	(G\otimes M)^{\otimes c} \cong G^{\otimes c}\otimes M^{\otimes c} \longrightarrow G\otimes M^{\otimes c}.
	\]
	Thus, for every $s\in H^0(Z,G\otimes M)$, the image of
	$s^{\otimes c}$ under this map is zero.
	
	Let $\eta$ be the generic point of $Z$. Since $\rho$ is
	étale over the locus of divisors with distinct points,
	the extension $k(Z)/k(\bP^n)$ is separable. Hence $G_\eta\cong
	k(X)\otimes_{k(\bP^n)}k(Z)$ is reduced; see \cite[\tagref{030U}]{stacks-project}.
	Choosing a generator $\ell$ of $M_\eta$ and writing
	$s_\eta=g\otimes\ell$, the preceding vanishing gives
	$g^c=0$. Thus $g=0$. Since $G\otimes M$ is locally free
	on the integral variety $Z$, it follows that $s=0$.
	Therefore
	\[
	H^0(X,B_q(-qH))=0.
	\]
	
	On the other hand, the unit inclusion
	$\cO_Z \longhookrightarrow G$ gives
	\[
	h^0(X,B_q) =h^0(Z,G\otimes L_q) \geq h^0(Z,L_q) =\prod_{\nu=1}^n(\nu q) =n!q^n>0.
	\]
	Hence $B_q$ is initialized with respect to $qH$.
	
	Finally, Using \eqref{eq320}, we get $H^0(X,E_e(-H))=0$ and $H^0(X,E_e)\neq0$. Thus $E_e$ is initialized with respect to $H$.
\end{proof}

\begin{remark} [The construction need not produce Ulrich bundles]
	\label{neednotul}
	The bundles $E_e$ constructed above need not be Ulrich.
	To see this, let $Q^3\subset\bP^4$ be a smooth quadric threefold,
	let $H$ be its hyperplane divisor, and let
	$\pi:Q^3\longrightarrow\bP^3$ be the projection from a point
	outside $Q^3$. Then $\pi$ is finite flat of degree $2$, and
	$\pi^*\cO_{\bP^3}(1)\cong\cO_{Q^3}(H)$.
	
	Choose coordinates in which the projection centre is
	$[0:0:0:0:1]$. Since this point is outside $Q^3$, the defining
	equation of $Q^3$ can be written as
	\[
	z^2+\ell(x_0,\ldots,x_3)z+g(x_0,\ldots,x_3)=0,
	\]
	where $\ell$ and $g$ are homogeneous of degrees $1$ and $2$,
	respectively. With $S=k[x_0,\ldots,x_3]$, the homogeneous
	coordinate ring is therefore free as a graded $S$-module, with
	basis $1,z$ of degrees $0,1$. Sheafifying gives
	\begin{equation}\label{e410}
		\pi_*\cO_{Q^3} \cong\cO_{\bP^3}\oplus\cO_{\bP^3}(-1).
	\end{equation}
	This argument is valid in every characteristic.
	
		Let $V$ be a vector bundle of positive rank on $\bP^3$, and let
	$q$ be a positive integer. By \eqref{e410} and the projection
	formula,
	\[
	H^i(Q^3,(\pi^*V)(-mqH))
	\cong H^i(\bP^3,V(-mq))\oplus H^i(\bP^3,V(-mq-1)).
	\]
	Thus, if $\pi^*V$ were Ulrich with respect to $qH$, both $V$
	and $V(-1)$ would be Ulrich with respect to
	$\cO_{\bP^3}(q)$. The Hilbert polynomial
	$P_V(t)=\chi(\bP^3,V(t))\in\mathbb Q[t]$, which has degree $3$ and leading
	coefficient $\operatorname{rk}(V)/6$, would then vanish at
	\[
	-q,\,-2q,\,-3q,\,-q-1,\,-2q-1,\,-3q-1.
	\]
	For any positive integer $q$, there are at least four distinct integers in this list which is impossible. Hence $\pi^*V$ is not Ulrich with respect to
	$qH$ for any positive integer $q$.
	
	In characteristic $p>0$, Proposition~\ref{prop34}, applied with
	$q=p^e$, shows that $F^e_*\pi^*V$ is not Ulrich with respect
	to $H$ for any $e\geq1$. In particular, the bundles $E_e$
	constructed above on $Q^3$ are not Ulrich.
		\end{remark}

\subsection{Relation with lim Ulrich sequences}
\label{weaklimrel}
The weakly Ulrich condition concerns an individual vector bundle and
requires exact vanishing in the ranges of Definition \ref{defweak}.
The lim Ulrich condition concerns a
sequence $(\mathcal E_e)$ of coherent sheaves. It requires
$h^0(\mathcal E_e)>0$, uniform vanishing for sufficiently negative
twists in degree zero and sufficiently positive twists in positive
degrees, and convergence of
$h^i(\mathcal E_e(tH))/h^0(\mathcal E_e)$ to zero throughout the
Ulrich vanishing range; see \cite[\S 6.6, p. 1239]{IMW22} and \cite[\S 2.4]{IMW23} . Such convergence does not force the
cohomology groups of any individual term to vanish. The following
examples show that a lim Ulrich sequence need not have weakly Ulrich
terms, and that a sequence of weakly Ulrich bundles need not be
lim Ulrich. Let $n\geq1$, and equip $\bP^n$ with $\cO_{\bP^n}(1)$.

\begin{example} [A lim Ulrich sequence with no weakly Ulrich terms]
	\label{limweakex1}
	Consider
	\[
	F_e=\cO_{\bP^n}^{\oplus e}\oplus\cO_{\bP^n}(-n-1), \qquad e\geq1.
	\]
	Then $h^0(F_e)=e$, and, for every fixed $i,t$,
	\[
	\frac{h^i(F_e(t))}{h^0(F_e)} =h^i(\cO_{\bP^n}(t)) +\frac{h^i(\cO_{\bP^n}(t-n-1))}{e}.
	\]
	Thus the normalized cohomology tables converge entrywise to the
	cohomology table of $\cO_{\bP^n}$. Moreover, for every $e\geq1$,
	\begin{align*}
		H^0(\bP^n,F_e(t))&=0 &&(t\leq-1),\\
		H^i(\bP^n,F_e(t))&=0 &&(i>0,\ t\geq1).
	\end{align*}
	These give the uniform vanishing bounds required in the definition,
	so $(F_e)$ is a lim Ulrich sequence. However,
	\[
	H^n(\bP^n,F_e)\cong k\ne0
	\]
	for every $e$, whereas \eqref{weak1}, with $i=n$ and $m=0$,
	requires this group to vanish. Hence none of the terms $F_e$ is
	weakly Ulrich \wrt $\cO_{\bP^n}(1)$.
\end{example}

\begin{example} [A sequence of weakly Ulrich bundles that is not lim Ulrich]
	\label{limweakex2}
	Consider
	\[
	G_e=\cO_{\bP^n}(1)^{\oplus e},\qquad e\geq1.
	\]
	The bundle $G_e$ has no intermediate cohomology, and
	\begin{align*}
		H^0(\bP^n,G_e(-m))&=0 &&(m\geq2),\\
		H^n(\bP^n,G_e(-m))&=0 &&(m\leq n-1).
	\end{align*}
	Thus every $G_e$ is weakly Ulrich. On the other hand,
	\[
	\frac{h^0(G_e(-1))}{h^0(G_e)} =\frac{e}{e(n+1)}=\frac1{n+1}.
	\]
	This ratio does not tend to zero, so $(G_e)$ is not a lim Ulrich
	sequence. 
\end{example}

\section{Existence of almost Ulrich bundles on threefolds}
In this section, we construct almost Ulrich bundles on smooth projective threefolds defined over an algebraically closed field in positive characteristic. The weakly Ulrich bundles $E_e$ constructed the previous section need not be almost Ulrich; see Remark \ref{nonalmostul}. 
Let $X$ be a smooth projective variety of dimension $n$ over an algebraically closed field $k$ of positive characteristic. Let $H$ be a very ample divisor. Write $E(m):=E \otimes \cO_X(mH)$.
\begin{definition}
A vector bundle $E$ is \emph{arithmetically Cohen-Macaulay}, or \emph{ACM}, with respect to $H$ if
$H^i(X,E(m))=0$ for every $m\in\bZ$ and every $0<i<n$.
\end{definition}

Almost Ulrich bundles \wrt a polarization $H$ are both ACM and weakly Ulrich.  Recall their definition from \cite{Jos21}.
\begin{definition} \label{almostul} \cite[Section 2, p. 3]{Jos21}
	A vector bundle $E$ on $X$ is called almost \ul \wrt $H$ if the following conditions are satisfied:
	\begin{itemize}
		\item [(1)] $H^n(X, E(-mH))=0$ for all integers $m \le (n-1)$.
		\item  [(2)] $H^j(X, E(mH))=0$ for $1 \le j \le n-1$ and for all $m \in \bZ$ (i.e., $E$ is ACM).
		\item  [(3)] $H^0(X, E(-mH))=0$ for all integers $m \ge 2$.
	\end{itemize}
\end{definition}

\begin{lemma}\label{almostullemcrit}
	Let $E$ be a vector bundle on $X$. Then $E$ is almost Ulrich \wrt $H$ if and only if $E$ is ACM \wrt $H$ and satisfies
	\[
	H^0(X,E(-2))=0,\qquad H^n(X,E(1-n))=0.
	\]
\end{lemma}

\begin{proof}
	Multiplication by a nonzero section of $\cO_X(1)$ injects $E(m-1)$ into
	$E(m)$. This induces an injection $H^0(X, E(m-1)) \longhookrightarrow H^0(X, E(m))$, and a surjection $H^n(X, E(m-1)) \twoheadrightarrow H^n(X, E(m))$. The assertion follows from these two considerations.
\end{proof}

\begin{remark} [The weakly Ulrich bundles $E_e$ need not be almost Ulrich]
	\label{nonalmostul}
	The weakly Ulrich bundle $E_e$ constructed in the previous section need not be almost Ulrich can be seen as follows. We consider the degree $2$ map $\pi: X=Q^3 \longrightarrow \bP^3$ from Remark \ref{neednotul}. For $Z=(\mathbb P^1)^3$, it follows that
	\[
	G=\rho^*\pi_*\mathcal O_X
	\cong
	\mathcal O_Z\oplus\mathcal O_Z(-1,-1,-1).
	\]
	For $n=3$, $L_q=\mathcal O_Z(q-1,2q-1,3q-1)$. The construction from the preceding section gives
	\[
	H^2(X,E_e(-2H)) \cong H^2\!\left(Z,\mathcal O_Z(-q-1,-1,q-1)\right)\oplus H^2\!\left(Z,\mathcal O_Z(-q-2,-2,q-2)\right).
	\]
	The first summand vanishes because its second factor is $\mathcal O_{\mathbb P^1}(-1)$. Since $q=p^e \ge2$, the K\"{u}nneth formula identifies the second summand with
	\[
	H^1(\mathbb P^1,\mathcal O(-q-2)) \otimes H^1(\mathbb P^1,\mathcal O(-2)) \otimes H^0(\mathbb P^1,\mathcal O(q-2)).
	\]
	Its dimension is $(q+1)\cdot1\cdot(q-1)=q^2-1$.  Therefore, $h^2(X,E_e(-2H))=q^2-1=p^{2e}-1>0$ for $e \ge 1$. Thus, on the quadric $Q^3$, the bundles $E_e$ are not ACM, and hence are not almost Ulrich.
\end{remark}

\subsection{Construction of an acyclic bundle} 
Throughout this subsection, let $X$ be a projective variety of dimension $n \ge 1$ over an algebraically closed field of arbitrary characteristic, and let $H$ be a very ample divisor on $X$. Write $E(m):=E \otimes \cO_X(mH)$. We construct acyclic vector bundles on $X$, without assuming that $X$ is smooth.

\begin{lemma} \label{lem1}
	For every sufficiently large integer $a$, and for every sufficiently large
	integer $b$ depending on $a$, the vector spaces
	\[
	V=H^0(X,\cO_X(a)),\quad W=H^0(X,\cO_X(b)),\quad
	T=H^0(X,\cO_X(b-a))
	\]
	have dimensions $r,s,t$ satisfying the following conditions:
	\begin{enumerate}
		\item $H^i(X,\cO_X(a))=H^i(X,\cO_X(b))=0$ for $i>0$;
		\item the multiplication map $\mu:V\otimes T\longrightarrow W$ is surjective;
		\item $r\ge2$, $rt\ge(r-1)(s+1)$, and $s-r\ge n$.
	\end{enumerate}
\end{lemma}

\begin{proof}
	Choose $a\ge1$ so that the first vanishing holds. Since $H$ is very ample,
	$\cO_X(a)$ is globally generated and $r\ge2$. Let $M_a$ be the kernel in
	the evaluation sequence
	\[
	0\longrightarrow M_a\longrightarrow V\otimes\cO_X \longrightarrow\cO_X(a)\longrightarrow0.
	\]
	After tensoring by $\cO_X(b-a)$, Serre vanishing gives
	$H^1(X,M_a(b-a))=0$ for $b\gg a$, and hence (ii).
	It also gives the second vanishing in (i).
	
	Let $P(m)=\chi(X,\cO_X(m))$. For $b$ sufficiently large,
	$s=P(b)$ and $t=P(b-a)$. If the leading term of $P$ is $\alpha m^n$,
	then $\alpha>0$ and
	\[
	rt-(r-1)(s+1)=rP(b-a)-(r-1)\bigl(P(b)+1\bigr)=\alpha b^n+O(b^{n-1}).
	\]
	This tends to $+\infty$, because $n\ge1$. Also $s-r\longrightarrow+\infty$.
	Thus all the required inequalities hold after increasing $b$.
\end{proof}

Fix $a,b$ as in Lemma \ref{lem1}. Let $\cA=M_{r\times s}(T)$ denote the vector space of $r \times s$ matrices with coefficients in $T=H^0(X,\cO_X(b-a))$. Then $\dim \cA=rst$. An element $M=(m_{ij})$ of $\cA$ defines a sheaf morphism
\begin{equation}\label{eq:phi}
	\varphi_M:\cO_X(a)^{\oplus s}\longrightarrow\cO_X(b)^{\oplus r}
\end{equation}
and a linear map
\[
\Phi_M=H^0(\varphi_M):V^{\oplus s}\longrightarrow W^{\oplus r}.
\]
The source and target of $\Phi_M$ both have dimension $rs$. We prove the following proposition using an incidence dimension count similar to that in the proof of \cite[Theorem 2.6]{coskun2025stability}.

\begin{proposition}\label{secisomlocus}
	The set $U_{\mathrm{sec}}=\{M\in\cA:\Phi_M\text{ is an isomorphism}\}$
	is a nonempty Zariski open subset of $\cA$.
\end{proposition}

\begin{proof}
	Note that $\Phi_M$ is an isomorphism if and only if $\det \Phi_M$ is nonzero, which is an open condition. Thus, $U_{\mathrm{sec}}$ is Zariski open in $\cA$.  To prove nonemptiness, consider the closed incidence variety
	\[
	I=\{(M,[z])\in\cA\times \bP(V^{\oplus s}):\Phi_M(z)=0\}.
	\]
	For a nonzero tuple $z=(z_1,\ldots,z_s)$, let
	$V_z=\langle z_1,\ldots,z_s\rangle\subseteq V$ and
	$W_z=\mu(V_z\otimes T)\subseteq W$.
	For fixed $z$, each row of $M$ induces the linear map
	\[
	T^{\oplus s}\longrightarrow W,\qquad
	(t_1,\ldots,t_s)\longmapsto\sum_j t_jz_j,
	\]
	whose image is $W_z$. Since the rows vary independently, the condition
	$\Phi_M(z)=0$ has codimension exactly $r\dim W_z$ in $\cA$.
	
	If $V_z=V$, then $W_z=W$ by Lemma \ref{lem1}(ii).
	The fibre of $I$ over $[z]$ has dimension $\dim \cA-rs$, and the locus of such
	$[z]$ has dimension at most $rs-1$. This part of $I$ therefore has
	dimension at most $\dim \cA-1$.
	
	Suppose now that $V_z\ne V$. Such tuples lie in $U^{\oplus s}$ for some
	hyperplane $U\subset V$. Hyperplanes in $V$ vary in dimension $r-1$, whereas
	$\bP(U^{\oplus s})$ has dimension $(r-1)s-1$. Hence the locus of these
	$[z]$ has dimension at most
	\[
	(r-1)(s+1)-1.
	\]
	At least one $z_j$ is nonzero. Since $X$ is integral, multiplication by
	this section injects $T$ into $W$: a nonzero section of an invertible
	sheaf is nonzero at the generic point and is a nonzerodivisor.
	Consequently $\dim W_z\ge t$. This part of $I$ has dimension at most
	\[
	\dim \cA-rt+(r-1)(s+1)-1\le \dim \cA-1.
	\]
	It follows that $\dim I\le \dim \cA-1$.
	
	The projection $I \longrightarrow \cA$ is proper, because its other factor
	is projective. Its image is thus a proper closed subset. Outside that image,
	$\Phi_M$ is injective and hence is an isomorphism.
\end{proof}

\begin{proposition} \label{surjlocus}
	The set $U_{\mathrm{surj}}=\{M\in\cA:\varphi_M\text{ is surjective on }X\}$ is a nonempty Zariski open subset of $\cA$.
\end{proposition}

\begin{proof}
	The line bundle $\cO_X(b-a)$ is globally generated. Therefore, for each
	closed point $x \in X$, evaluation of matrices gives a surjection 
	$M_{r\times s}(T) \twoheadrightarrow M_{r\times s}(k)$ after choosing trivializations of the
	fibres of $\cO_X(b-a)$ at $x$. The locus of matrices of rank less than $r$ has codimension
	$s-r+1$. Indeed, the locus of matrices of rank at most $r-1$ has dimension
	$(r-1)(s+1)$ in the $rs$-dimensional space of matrices.
	
	The closed incidence locus
	\[
	J=\{(M,x)\in\cA\times X:\rk(\varphi_M\otimes k(x))<r\}
	\]
	therefore has dimension at most $\dim \cA+n-(s-r+1)\le \dim \cA-1$.
For the first projection, the image of $J$ in $\cA$ is closed because $X$ is projective, and is proper by
the dimension bound. The complement is precisely $U_{\mathrm{surj}}$.
\end{proof}

\begin{theorem}\label{acyclic}
	Let $X$ be a projective variety of dimension $n\ge1$ over an
	algebraically closed field of arbitrary characteristic. Let $H$ be a
	very ample divisor. For $a,b$ as in Lemma \ref{lem1}, a general matrix
	$M\in\cA$ gives an exact sequence
	\begin{equation} \label{kerbun}
		0\longrightarrow B_M\longrightarrow\cO_X(a)^{\oplus s}
		\xrightarrow{\varphi_M}\cO_X(b)^{\oplus r}\longrightarrow0,
	\end{equation}
	where $B_M$ is a vector bundle of positive rank $s-r$ and
	$H^\bullet(X,B_M)=0$.
\end{theorem}

\begin{proof}
	The affine space $\cA$ is irreducible. By
	Propositions \ref{secisomlocus} and \ref{surjlocus}, the intersection
	$U_{\mathrm{sec}}\cap U_{\mathrm{surj}}$ is a nonempty open subset.
	It has a $k$-point because $k$ is algebraically closed.
	For every matrix in this intersection, the kernel in \eqref{kerbun}
	is locally free: locally, the surjection onto a free module splits.
	Its rank is $s-r\ge n>0$.
	The map on global sections is an isomorphism, and both other terms in
	\eqref{kerbun} have zero higher cohomology by Lemma \ref{lem1}(i).
	The long exact cohomology sequence gives the assertion.
\end{proof}

\subsection{Construction of almost Ulrich bundles on threefolds}
Let $X$ be a smooth projective variety of dimension $n \ge 1$ defined over an algebraically closed field of positive characteristic. Let $H$ be a very ample divisor. Fix a nonzero acyclic vector bundle $B$ on $X$ from Theorem \ref{acyclic}.

\begin{lemma}\label{serrevangg}
	There is an integer $N\ge1$ such that, for every integer $\ell\ge N$,
	\[
	\begin{aligned}
		&H^i(X,B(\ell))=0\quad(i>0), &&h^0(X,B(\ell))>0,\\
		&H^i(X,B(-\ell))=0\quad(i<n), &&h^n(X,B(-\ell))>0.
	\end{aligned}
	\]
\end{lemma}

\begin{proof}
		By Serre vanishing and eventual global generation \cite[III, Theorem 5.2 and II, Theorem 5.17]{hartshorne1977algebraic}, applied to $B$ and
		$B^\vee\otimes\omega_X$, both sheaves have vanishing higher cohomology
		and nonzero global sections after sufficiently positive twists. Serre
		duality gives
		\[
		H^i(X,B(-\ell))^\vee \cong H^{n-i}(X,B^\vee\otimes\omega_X(\ell)),
		\]
		which yields the remaining assertions for $\ell\gg0$.
\end{proof}

Let $F:X \longrightarrow X$ denote the absolute Frobenius morphism, and let $F^e:X \longrightarrow X$ denote its $e$-fold iterate. Since $X$ is smooth, $F^e$ is finite flat of degree $p^{en}$. Let $q=p^e$. Then $C_e=F^e_*B$ is a vector bundle on $X$ of rank $q^n \rk(B)$.

\begin{proposition} \label{cohoFdirect}
Let $H$ be a very ample divisor on $X$. For sufficiently large $e$, for every $m \in \bZ$ and every $i \ge 0$,
	\begin{equation} \label{Fdirecte1}
		H^i(X,C_e(m))\ne0 \quad\Longleftrightarrow\quad \bigl(m>0\text{ and }i=0\bigr)
		\ \text{or}\  \bigl(m<0\text{ and }i=n\bigr).
	\end{equation}
	In particular, $C_e$ is ACM \wrt $H$, and $C_e(m)$ is acyclic if and only if $m=0$.
\end{proposition}

\begin{proof}
	The projection formula and finiteness give, for every $m$ and $i$,
	\begin{equation}\label{Fdirecte2}
		C_e(m)\cong F^e_*(B \otimes (F^e)^* \cO_X(mH)) \cong F^e_*\bigl(B(mq)\bigr),\qquad
		H^i(X,C_e(m))\cong H^i(X,B(mq)).
	\end{equation}

	Choose $e$ such that $q \ge N$, with $N$ as in Lemma \ref{serrevangg}. If $m>0$ then
	$mq\ge N$, and if $m<0$ then $mq \le -N$. The two assertions of that lemma
	apply uniformly to every nonzero integer $m$. For $m=0$ use $H^\bullet(X,B)=0$. This proves \eqref{Fdirecte1}.
\end{proof}

\begin{theorem} \label{almostulthm}
	Let $X$ be a smooth projective variety over an algebraically closed field of positive characteristic. For each $c \in \bZ$, set $E_{e,c}=C_e(c)$. Then for all sufficiently large $e$, the following hold simultaneously for every integer $c$:
\begin{enumerate}
		\item $E_{e,c}$ is ACM;
		\item $E_{e,c}$ is almost Ulrich if and only if $n-1\le c\le2$;
		\item $E_{e,c}$ is Ulrich if and only if $n=1$ and $c=1$;
		\item $E_{e,c}$ is initialized if and only if $c=1$.
	\end{enumerate}
\end{theorem}

\begin{proof}
	The ACM assertion follows from Proposition \ref{cohoFdirect}.
	The same proposition gives
	\[
	\begin{aligned}
		H^0(X,E_{e,c}(-2))=0&\quad\Longleftrightarrow\quad c\le2,\\
		H^n(X,E_{e,c}(1-n))=0&\quad\Longleftrightarrow\quad c\ge n-1.
	\end{aligned}
	\]
	Lemma \ref{almostullemcrit} proves (ii).
	
	For (iii), the Ulrich condition requires
	$C_e(c-j)$ to be acyclic for every $j=1,\ldots,n$.
	By Proposition \ref{cohoFdirect}, this is equivalent to
	$c-j=0$ for every such $j$. This holds exactly when $n=1$ and $c=1$.
	
	Finally, \eqref{Fdirecte1} gives
	$h^0(X,E_{e,c})>0$ exactly when $c\ge1$, and
	$h^0(X,E_{e,c}(-1))=0$ exactly when $c\le1$.
	Their conjunction is $c=1$.
\end{proof}

\begin{proof}[Proof of Theorem \ref{thm4}]
	When $X$ is a smooth projective threefold, substituting $n=3$ in Theorem \ref{almostulthm} shows that $E_{e, 2}=C_e(2)=(F^e_*B)(2)$ is almost Ulrich \wrt $H$ for all sufficiently large $e$.
\end{proof}

\subsection{Reduction from algebraically closed field to arbitrary field for existence}
\label{extensionarbifields}
The existence results for Ulrich, weakly Ulrich, and almost Ulrich bundles in Theorems \ref{thm1}, \ref{thm2}, \ref{thm3}, and \ref{thm4} extend to arbitrary fields under the respective characteristic assumptions. This follows from a standard descent argument that works in arbitrary characteristic.

Let $k_0$ be a field, and let $\overline{k_0}$ be an algebraic closure of $k_0$. Let $X$ be a smooth projective variety over $k_0$, and let $H$ be a very ample divisor on $X$. For every field extension $K/k_0$, write $X_K:=X \times_{\spec k_0} \spec K$,
and let $H_K$ denote the pullback of $H$ to $X_K$. We show that if $(X_{\overline{k_0}},H_{\overline{k_0}})$ admits an Ulrich, weakly Ulrich, or almost Ulrich bundle, then $(X,H)$ admits a bundle with the same cohomological property, possibly of larger rank.

Let $\overline{\mathcal E}$ be such a bundle on $X_{\overline{k_0}}$. Since $\overline{k_0}$ is the filtered union of its finite subextensions over $k_0$ and $\overline{\mathcal E}$ is finitely presented, there exist a finite extension $L/k_0$ and a finitely presented sheaf $\mathcal E_L$ on $X_L$ such that $g^*\mathcal E_L\simeq\overline{\mathcal E}$,
where $g:X_{\overline{k_0}}\longrightarrow X_L$ is the natural projection; see \cite[\tagref{01ZR}]{stacks-project}. Since $g$ is faithfully flat and affine, local freeness descends along $g$. Thus $\mathcal E_L$ is a vector bundle by \cite[\tagref{05B2}]{stacks-project}.

For every $i\ge0$ and every $t\in\mathbb Z$, flat base change gives
\begin{equation}\label{bce1}
	H^i \left(X_L, \mathcal E_L(tH_L)\right) \otimes_L\overline{k_0} \simeq
	H^i \left(X_{\overline{k_0}}, \overline{\mathcal E}(tH_{\overline{k_0}})\right).
\end{equation}
By faithful flatness, $H^i(X_L,\mathcal E_L(tH_L))=0$ if and only if  $H^i(X_{\overline{k_0}},\overline{\mathcal E}(tH_{\overline{k_0}}))=0$.
Consequently, $\mathcal E_L$ is Ulrich, weakly Ulrich, or almost Ulrich \wrt $H_L$, respectively. 

Now consider the natural projection $p: X_L\longrightarrow X$. It is finite flat of degree $[L:k_0]$, so $\mathcal E:=p_* \mathcal E_L$
	is a vector bundle on $X$. The projection formula gives $\mathcal E(tH)\simeq p_* \bigl(\mathcal E_L(tH_L)\bigr)$.
	Since $p$ is finite, we obtain
	\[
	H^i(X,\mathcal E(tH)) \simeq H^i(X_L,\mathcal E_L(tH_L))
	\]
	for every $i\ge0$ and every $t\in\mathbb Z$. Thus all the cohomology vanishings defining the Ulrich, weakly Ulrich, or almost Ulrich property are preserved, and $\mathcal E$ has the corresponding property \wrt $H$. If $\mathcal E_L$ has constant rank $r$, then $\rk(\mathcal E)=r[L:k_0]$.

\subsection*{Acknowledgement}
The author is grateful to Adrian Langer for answering his questions about moduli spaces of sheaves in positive characteristic. The author is supported by a postdoctoral fellowship from IISER Pune.

\bibliographystyle{plainurl}
\bibliography{ulrich.bib}

	\end{document}